\documentclass[MSNbibl,number,citesort,dvips]{arxbj}
\usepackage{upgreek}
\usepackage{mathrsfs}

\aid{0}
\volume{20}
\issue{2}
\pubyear{2014}
\firstpage{395}
\lastpage{415}
\doi{10.3150/12-BEJ491} %kopijuoti is PTS

\makeatletter

\newcommand{\rrvert}{\vert}
\newcommand{\llvert}{\vert}
\newcommand{\E}{\mathbb{E}}
\newcommand{\pr}{\mathbb P}
\def\s{^\star}
\newcommand{\toi}{\to\infty}

\newcommand{\refs}[1]{(\ref{#1})}
\def\R{\mathbb{R}}
\newtheorem{thmm}{Theorem}
\newtheorem{lem}{Lemma}
\newtheorem{cor}{Corollary}
\newtheorem{prop}{Proposition}
\newremark{rem}{Remark}

\renewcommand{\epsilon}{\varepsilon}
\makeatother

\begin{document}
\begin{frontmatter}

\title{On the tail asymptotics of the area swept under the Brownian storage graph}
\runtitle{On the tail asymptotics of the area swept under the Brownian storage graph}

\begin{aug}
%%%% inicialai - be tarpu
\author[a]{\fnms{Marek} \snm{Arendarczyk}\thanksref{a,e1}\ead[label=e1,mark]{marendar@math.uni.wroc.pl}},
\author[a]{\fnms{Krzysztof} \snm{D\c{e}bicki}\corref{}\thanksref{a,e2}\ead[label=e2,mark]{debicki@math.uni.wroc.pl}}
\and
\author[b]{\fnms{Michel} \snm{Mandjes}\thanksref{b}\ead[label=e3]{M.R.H.Mandjes@uva.nl}}
\runauthor{M. Arendarczyk, K. D\c{e}bicki and M. Mandjes} %% auto

\address[a]{Mathematical Institute, University of
Wroc\l aw, pl. Grunwaldzki 2/4, 50-384 Wroc\l aw, Poland.\\
\printead{e1,e2}}
\address[b]{Korteweg-de Vries Institute for Mathematics,
University of Amsterdam, The Netherlands; Eurandom,
Eindhoven University of Technology, The Netherlands; CWI,
Amsterdam, The Netherlands.
\printead{e3}}

\end{aug}

% HISTORY:
\received{\smonth{6} \syear{2011}}
\revised{\smonth{8} \syear{2012}}

% ABSTRACT
%
\begin{abstract}In this paper, the area swept under the workload graph
is analyzed: with
{$\{Q(t)\dvt t\ge0\}$ denoting the stationary workload
process}, the asymptotic behavior of
\[
\pi_{T(u)}(u):={\pr} \biggl(\int_0^{T(u)}
Q(r)\,\mathrm{d}r>u \biggr)
\]
{is} analyzed. Focusing on regulated Brownian motion, first the
exact asymptotics of $\pi_{T(u)}(u)$ are given for the case that
$T(u)$ grows slower than $\sqrt{u}$, and then logarithmic
asymptotics for (i) $T(u)= T\sqrt{u}$ (relying on sample-path
large deviations), and (ii) $\sqrt{u}=\mathrm{o}(T(u))$ but $T(u)=\mathrm{o}(u)$.
Finally, the Laplace transform of the residual busy period are
given in terms of the Airy function.
%%, leading to explicit expressions for the mean and variance.
%the paper is devoted to regulated \textit{fractional} Brownian
%motion, which is intrinsically harder; for the short time-scale
%regime the exact asymptotics of $\pi_{T(u)}(u)$ are presented.\fi
\end{abstract}

% KEYWORDS
% visi is mazosios raides ir pagal abecele
%

\end{frontmatter}

%s1 #&#
\section{Introduction}
Queueing models form an important branch within applied probability,
having applications in production, storage, and inventory systems, as
well as in communication networks. At the same time, there is a strong
link with various models that play a crucial role in finance and risk
theory, see, for instance, \cite{Kyprianou06}.

In more formal terms, the workload process of a queue is commonly
defined as follows. Let $(X(t))_{t\in{\mathbb R}}$ be a stochastic
process, that is often assumed to have stationary increments; without
loss of generality we assume it has zero mean. Let $c>0$ be the drain
rate of the queue. Then the corresponding workload process $(Q(t))_{t\in
{\mathbb R}}$ is defined through
\[
Q(t) = \sup_{s\le t} X(t) - X(s) - c(t-s).
\]
A sizable body of literature is devoted to the analysis of the
probabilistic properties of this workload process, both in terms of its
stationary behavior and its transient characteristics.

One of the key metrics of the queueing system under consideration is
the mean stationary workload. In many situations, this cannot be
computed explicitly, and one then often resorts to simulation. A
commonly used estimator is
\[
\bar Q_T := \frac{1}{T} \sum_{i=1}^T
Q(i);
\]
one could set up the situation such that at time 0 the queue has
already run for a substantial amount of time, such that one can safely
assume the workload is in stationarity. In the simulation literature,
this type of estimators (and related ones) have been analyzed in
detail; see, for example,~\cite{Asmussen07}. Results are in terms of laws of
large numbers and central limit theorems.

Recently, attention shifted to the large deviation properties of the
above type of estimators. It is observed that the subsequent
observations are in general dependent, which considerably complicates
the analysis. More specifically, standard large-deviations techniques
do not apply here; the G\"artner--Ellis theorem \cite{Dembo98}, that
allows only a mild dependence between the increments, is therefore not
of any use.
Even in cases in which the correlation of the stationary workload
exhibits roughly exponential decay (being a manifestation of the
queue's input process having short-range dependent properties), it
turns out that the probability of the sample mean $\bar Q_T$ deviating
from the mean stationary workload, say $q$, under quite general
circumstances, does \textit{not} decay exponentially.

Let us consider a few more detailed results.
In a random walk setting (i.e., in which $Q(0)=0$ and $Q(t+1) = \max\{
Q(t)+Y(t),0\}$ for an i.i.d. sequence $Y(t)$), Meyn \cite{Meyn06,Meyn08} proved an intriguing (asymmetric) result. `Below the mean'
there is, under mild regularity assumptions, exponential decay, in that
\[
\limsup_{T\toi}\frac{1}{T}\log\pr(\bar Q_T\le a) <0
\]
for each $a<q$, whereas `above the mean' there is `subexponential
decay', that is,
\[
\lim_{T\toi}\frac{1}{T}\log\pr(\bar Q_T\ge a) = 0
\]
for each $a>q$. Subsequently, Duffy and Meyn \cite{Duffy10}, proved
that the right scaling was \textit{quadratic}, in the sense that in
their setting ${T}^{-2} \sum_{i=1}^T Q(i)$ satisfies a large deviations
principle with a nontrivial rate function. The square can intuitively
be understood from the fact that one essentially considers the right
scaling for the \textit{area} under the graph of the workload.

The above motivates the interest in tail probabilities of the type
\[
\pi_{T(u)}(u) := \pr \biggl(\int_0^{T(u)}
Q(t)\,\mathrm{d}t > u \biggr)
\]
for various types of interval lengths $T(u)$, and $u\toi$;
here the workload is assumed to be in stationarity at time 0.
As indicated above, for $T(u)$ be in the order of $\sqrt{u}$ and the
queue's input process having i.i.d. increments, the tail probability
$\pi_{T(u)}(u)$ decaying roughly like $\exp(-\alpha\sqrt{u})$ for some
$\alpha\in(0,\infty)$. On the other hand, for the case $u=\mathrm{o}(T(u))$ it
is seen that $\pi(u)$ tends to 1 for $u$ large.

The queueing system we consider in this paper is \textit{reflected}
(or: regulated) \textit{Brownian motion}, also referred to as \textit
{Brownian storage}; this means that the driving process $(X(t))_{t\in
{\mathbb R}}$
is {a} (standard) Brownian motion. In more detail, our contributions
are the following.
\begin{itemize}
\item We first, in Section \ref{BMexactas}, consider the
\textit{short timescale} regime, that is, we assume $T(u)=\mathrm{o}(\sqrt{u})$.
The main intuition here is that, in this regime, with overwhelming
probability the queue does not idle in $[0,T(u)]$, and as a
consequence, $Q(s)$ behaves as $Q(0)+X(s)-cs$ for $s\in[0,T(u)].$ This
essentially enables us to compute the so-called \textit{exact
asymptotics} of $\pi_{T(u)}(u)$, that is, we find an explicit function
$\varphi(u)$ such that $\pi_{T(u)}(u)/\varphi(u)\to1$ as $u\toi.$
\item The second contribution concerns the \textit
{intermediate timescale} regime, in which $T(u)$ is proportional to
$\sqrt{u}.$ As a function of this proportionality constant, we
determine in Section \ref{BMlogas} the decay rate
\[
-\alpha=\lim_{u\to\infty}\frac{1}{\sqrt{u}} \log\pi_{T(u)}(u),
\]
such that $\pi_{T(u)}(u)$ roughly looks like $\exp(-\alpha\sqrt{u})$
for $u$ large.
A crucial observation is that the probability under study can be
translated into a related probability in the so-called many-sources
regime. This means that sample-path large deviations for Brownian
motion can be applied here, for example, Schilder's theorem. Apart from
determining the decay rate, also the associated most likely path is
identified, complementing results in \cite{Duffy10}.
\item Section \ref{BMlong} considers the \textit{long
timescale}, that is $\sqrt{u}=\mathrm{o}(T(u))$ but $T(u)=\mathrm{o}(u)$. Relying on the
intuition that essentially one `big' busy period causes the rare event
under consideration, we prove that (like in the intermediate timescale
regime) $\pi_{T(u)}(u)$ roughly decays like $\exp(-\alpha\sqrt{u})$
for some constant $\alpha>0$. The proof techniques are reminiscent of
those used to establish an analogous property in the M$/$M$/$1 queue~\cite
{BlanchetGlynnMeyn11}.
\item We then consider in Section \ref{BMrunBP} the integral
over the \textit{remaining busy period}
(rather than a given horizon $T(u)$), again with Brownian motion input
(cf. the results for `traditional' single-server queues in \cite{Borovkov03a}).
It turns out to be possible to explicitly compute its Laplace
transform, in terms of the so-called Airy function,
which also enables closed-form expressions for the corresponding mean value.
%$t^{-1} \int_0^t Q(s)\mathrm{d}s$, in the spirit of the results for the
%M/G/1 queue in \cite{Borovkov03a}.
\end{itemize}

%s2 #&#
\section{Notation and model description}\label{snotation}

Let the stochastic process $\{B(t)\dvt t\in\R\}$ be a standard Brownian motion
{(i.e., ${\mathbb E}B(t)=0$ and ${\operatorname{\mathbb{V}ar}} B(t) = t$);
${\mathscr N}$ denotes a standard Normal random variable.}

In this paper, we consider a fluid queue fed by $B(\cdot)$ and drained
with a constant rate $c>0$.
Let $\{Q(t)\dvt t\in{\mathbb R}\}$ denote the stationary buffer content
process, that is, the unique stationary solution of the following
\textit{Skorokhod problem}:
\begin{longlist}
\item[{S1}] $Q(t) = Q(0)+B(t)-ct + L(t),$ for $ t\geq0$;
\item[{S2}] $Q(t)\geq0,$ for $t\geq0$;
\item[{S3}] $L(0)=0$ and $L$ is nondecreasing;
\item[{S4}] $\int_0^\infty Q(s) \,\mathrm{d}L(s) = 0$.
\end{longlist}

We recall that the solution to the above Skorokhod problem is
\[
Q(t)=\sup_{s\le t}\bigl(B(t)-B(s)-c(t-s)\bigr).
\]
The primary focus of this paper concerns the tail asymptotics
\[
\pi_{T(u)}(u):={\mathbb P} \biggl(\int_0^{T(u)}
Q(r)\,\mathrm{d}r>u \biggr)
\]
for functions $T(\cdot)\dvtx {\mathbb R}\to{\mathbb R}_+$.%; here the
%workload process is assumed to be in stationarity at time 0.
%In this paper we use the following additional notation.
%%
%%
%For given $H\in(0,1]$, write
%$\mathscr{H}_H(S):=\E\exp (\sup_{t \in[0,S]}
%By $\mathscr{H}_H$ we denote the \emph{Pickands's constant}:
%%
%%
%It is standard that, with $\mathscr{N}$ denoting the
%standard normal random variable,
%%
%%e1 #&#
%
%%

%s3 #&#
\section{Short timescale}\label{BMexactas}
In this section, we focus on the analysis of $\pi_{T(u)}(u)$ {as}
$u\to\infty$ and $T(u)=\mathrm{o}(\sqrt{u})$. The main intuition in this
timescale is that with overwhelming probability the queue does not
idle in $[0,T(u)]$. Therefore, $Q(r)$ essentially behaves as
$Q(0)+B(r)-cr$ for $r\in[0,T(u)]$, so that $\pi_{T(u)}(u)$ looks
like ($u$ large)
\[
\pr \biggl(\int_0^{T(u)} \bigl[Q(0) + B(r) - cr
\bigr] \,\mathrm{d}r > u \biggr).
\]
This idea is formalized in the following theorem.

%th1 #&#
\begin{thmm} \label{thbm}
Let $T(u) = \mathrm{o}(\sqrt{u})$. Then, as $u\toi$,
\[
\pi_{T(u)} (u) = \exp \biggl( -\frac{2cu}{T(u)} - \frac{1}{3}
c^2 T(u) \biggr) \bigl(1 + \mathrm{o}(1)\bigr).
\]
\end{thmm}

The following lemma plays an important role in the proof of Theorem \ref{thbm}.
%
%le1 #&#
\begin{lem} \label{lembm}
For any $T(\cdot)\dvtx {\mathbb R}\to{\mathbb R}_+$, as $u\toi$,
\[
\pr \biggl(\int_0^{T(u)} \bigl[Q(0) + B(r) - cr
\bigr] \,\mathrm{d}r > u \biggr) = \exp \biggl( -\frac{2cu}{T(u)} -
\frac{1}{3} c^2 T(u) \biggr) \bigl(1 + \mathrm{o}(1)\bigr).
\]
\end{lem}

\begin{pf}
Recalling that we assumed that the workload process is in steady-state
at time 0, it is well-known that
%
%e2 #&#
\begin{equation}
\label{distrQ} \pr\bigl(Q(0)>u\bigr)=\exp(-2cu),
\end{equation}
see, for example, Section 5.3 in \cite{Mandjes07}.
The distributional equality, for $T(u)>0$,
%
%e3 #&#
\begin{equation}
\label{intN} \int_0^{T(u)} B(t) \,\mathrm{d}t
\stackrel{\mathrm{d}} {=} \biggl(\frac{T(u)}{3} \biggr)^{1/2}
\mathscr{N}
\end{equation}
implies
\[
\pr \biggl(\int_0^{T(u)} \bigl[Q(0) + B(r) - cr
\bigr] \,\mathrm{d}r > u \biggr) = \pr \biggl(T(u) Q(0) + \frac{T(u)^{3/2}}{\sqrt{3}}
\mathscr{N} > u + \frac{1}{2} c T(u)^2 \biggr).
\]
Denote
\[
A_1(u) := \frac{ \sqrt{3}(u + ({1}/{2})c(T(u))^2 ) }{ (T(u))^{3/2}}.
\]

Integrating with respect to the distribution of $\mathscr{N}$, and
using (\ref{distrQ}), we obtain that
\begin{eqnarray*}
&&\pr \biggl(T(u) Q(0) + \frac{(T(u))^{3/2}}{\sqrt{3}} \mathscr{N} > u +
\frac{1}{2} c \bigl(T(u)\bigr)^2 \biggr)
\\
&&\quad= \frac{1}{\sqrt{2\uppi}} \int_{-\infty}^\infty
\pr \biggl( Q(0) > \frac{u}{T(u)} + \frac{1}{2} c T(u) - \biggl(
\frac{T(u)}{3} \biggr)^{1/2} x \biggr) \mathrm{e}^{-{x^2}/{2}}
\,\mathrm{d}x =I_1+I_2,
\end{eqnarray*}
with
\begin{eqnarray*}
I_1&:=& \frac{1}{\sqrt{2\uppi}} \exp \biggl(-\frac{2cu}{T(u)} -
c^2 T(u) \biggr) \int_{-\infty}^{A_1(u)} \exp
\biggl(-\frac{x^2}{2} - 2c \biggl( \frac{T(u)}{3} \biggr)^{{1} /{2}} x \biggr) \,\mathrm{d}x;
\\
I_2&:=& \frac{1}{\sqrt{2\uppi}} \int_{A_1(u)}^\infty
\exp \biggl( -\frac{x^2}{2} \biggr) \,\mathrm{d}x.
\end{eqnarray*}

\textit{Integral $I_1$}: First, rewrite
\[
I_1
= \frac{1}{\sqrt{2\uppi}} \exp \biggl(-\frac{2cu}{T(u)}
- \frac{1}{3} c^2 T(u) \biggr) \int_{-\infty}^{A_1(u)}
\exp \biggl(- \biggl(\frac{x}{\sqrt{2}} + A_2(u)
\biggr)^2 \biggr) \,\mathrm{d}x,
\]
where
$A_2(u) = c\sqrt{{2T(u)}/{3}}.$ Using the substitution
$
y := x + A_2(u)
$,
we obtain
\begin{eqnarray*}
I_1&=&\frac{1}{\sqrt{2\uppi}} \exp \biggl(-\frac{2cu}{T(u)} -
\frac{1}{3} c^2 T(u) \biggr) \int_{-\infty}^{A_1(u) + A_2(u)}
\exp \biggl(- \frac{y^2}{2} \biggr) \,\mathrm{d}y
\\
&=& \exp \biggl( -\frac{2cu}{T(u)} - \frac{1}{3} c^2 T(u)
\biggr) \bigl(1 + \mathrm{o}(1)\bigr)
\end{eqnarray*}
as $u \to\infty$.\vadjust{\goodbreak}

\textit{Integral $I_2$}:
\begin{eqnarray*}
I_2
&=& \frac{1}{\sqrt{2\uppi}A_1(u)} \exp \biggl( -
\frac{(A_1(u))^2}{2} \biggr) \bigl(1 + \mathrm{o}(1)\bigr)
\\
&=& \frac{ (T(u))^{3/2} }{ \sqrt{6\uppi}(u + ({1}/{2})c(T(u))^2 ) } \exp \biggl( -\frac{ 3\sqrt{6\uppi}(u + ({1}/{2})c(T(u))^2 ) }{
2(T(u))^3 }
\biggr)^2 \bigl(1 + \mathrm{o}(1)\bigr)
\\
&=& \mathrm{o} \biggl( \exp \biggl( -\frac{2cu}{T(u)} - \frac{1}{3}
c^2 T(u) \biggr) \biggr)
\end{eqnarray*}
as $u \to\infty$, where we used that
$\pr ( \mathcal{N}>x )\sim\frac{1}{\sqrt{2\uppi x}}\exp(-x^2/2)$
as $x\to\infty$.

This completes the proof.
\end{pf}

\begin{pf*}{Proof of Theorem \ref{thbm}} We establish upper and lower
bound separately.

\textit{Upper bound}:
%where
% L(t) = 0 \vee\sup_{s \in[0,t]} (-Q(0) - B(s) + cs).
We distinguish between the case that the queue has idled before $T(u)$,
and the case the buffer has been nonnegative all the time. We thus obtain
$\pi_{T(u)}(u)=P_1(u)+P_2(u),$
where
\begin{eqnarray*}
P_1(u)&:=&\pr \biggl( \int_{0}^{T(u)}
Q(r) \,\mathrm{d}r > u, L\bigl(T(u)\bigr) = 0 \biggr),
\\
P_2(u)&:=& \pr \biggl( \int_{0}^{T(u)}
Q(r) \,\mathrm{d}r > u, L\bigl(T(u)\bigr) > 0 \biggr).
\end{eqnarray*}
Due to {S1} and Lemma \ref{lembm}, as $u \to\infty$,
%
%e4 #&#
\begin{eqnarray}\label{rhs}
P_1(u) &\le& \pr \biggl(\int
_0^{T(u)} \bigl[Q(0) + B(r) - cr\bigr] \,\mathrm{d}r >
u \biggr)
\nonumber
\\[-8pt]
\\[-8pt]
\nonumber
%&=&
% \pr\left( \int_{0}^{T(u)} [Q(0) + B(t) - ct] dt > u \right)\\
 &=& \exp \biggl( -
\frac{2cu}{T(u)} - \frac{1}{3} c^2 T(u) \biggr) \bigl(1 + \mathrm{o}(1)
\bigr).
\end{eqnarray}
Moreover, for any $T(u)>0$,
\[
P_2(u) \le \pr \biggl( \sup_{s,t \in[0,T(u)]} \bigl[B(t) - B(s) - c(t -
s)\bigr] > \frac
{u}{T(u)} \biggr),
\]
realizing that for some epoch in $[0,T(u)]$ the workload has exceed
level $u/T(u)$, whereas for another epoch
it has been $0$.
According to the Borell inequality \cite{Adler90}, Theorem 2.1, in
conjunction with the self-similarity of Brownian motion, $P_2(u)$ is
majorized by
\begin{eqnarray*}
&& 2\exp \biggl(- \frac{  (({u}/{T(u)}) - \E [\sup_{s,t \in[0,T(u)]}
B(t) - B(s) - c(t - s) ] )^2}{2T(u)} \biggr)
\\
% \nonumber
% &\le&
% 2\exp
% \left(- \frac{ \left(\frac{u}{T(u)} - \E\left[\sup_{s,t \in[0,T(u)]}
%B(t) - B(s)\right]\right)^2}{2T(u)} \right) \\
&&\quad \le 2
\exp \biggl(- \frac{  (({u}/{T(u)}) - c T(u) - \sqrt{T(u)} \E [\sup_{s,t \in[0,1]} B(t) - B(s) ] )^2}{2T(u)} \biggr),
\end{eqnarray*}
which is negligible with respect to (\ref{rhs}) as $u \to\infty$. This
completes the proof of the upper bound.

\textit{Lower bound}: In view of
\[
\pr \biggl( \int_{0}^{T(u)} Q(t) \,\mathrm{d}t > u
\biggr)\ge\pr \biggl(\int_0^{T(u)} \bigl[Q(0) + B(r)
- cr\bigr] \,\mathrm{d}r > u \biggr),
\]
due to Lemma \ref{lembm} the proof is complete.
\end{pf*}

%s4 #&#
\section{Intermediate timescale}\label{BMlogas}
In this section, we consider the case of $T(u)$ being proportional to
$\sqrt{u}$: we set $T(u) = T\sqrt{u}$ {for some $T>0$}. The main result
of this section is given in the
following theorem, that describes the asymptotics of the probability
that the area until time $T\sqrt{u}$
exceeds $Mu$. It uses the following notation:
\[
\varphi(T,M):=\cases{
\frac{2}{3}\sqrt{6}
c\sqrt{cM}, &\quad$ \mbox{if $\sqrt{6M/c}<T$};$
\vspace*{2pt}\cr
2cM/T + c^2T/3,  &\quad $\mbox{else.}$}
\]
In this regime the intuition is that, in order to build up an area of
at least $u$, for relatively small values of $T$ the queue does not
idle with high probability, leading to an expression for the decay rate
that involves both $M$ and $T$. If, on the contrary, $T$ is somewhat
larger, then the most likely path is such that the queue starts off
essentially empty at time 0, to return to 0 before $T\sqrt{u}$, thus
yielding a decay rate that just depends on $M$.

%th2 #&#
\begin{thmm} \label{thinterm}
For all $T,M>0$, it holds that
%
%e5 #&#
\begin{equation}
\label{stated} - \lim_{u\to\infty}\frac{1}{\sqrt{u}}\log \pr \biggl(\int
_0^{T\sqrt{u}} Q(r)\,\mathrm{d}r \ge Mu \biggr) =
\varphi(T,M).
\end{equation}
\end{thmm}

We first observe that the probability under study can be translated
into a related probability in the so-called \textit{many-sources
regime}, as will be shown in Lemma \ref{lmscaling}.
Let $B^{(i)}(\cdot)$ be a sequence of independent standard Brownian motions.
Define
\[
\overline{B^{(n)}}(t) := \frac{1}{n} \sum
_{i=1}^n B^{(i)}(t),\qquad  Q^{(n)}(t) :=
\sup_{s\le t} \bigl(\overline{B^{(n)}}(t)- \overline
{B^{(n)}}(s)-c(t-s) \bigr).
\]

%le2 #&#
\begin{lem}\label{lmscaling}
For each $T,M > 0, n\in\mathbb{N}$
%it holds that
%
%e6 #&#
\begin{equation}
\label{eqlem2} \pr \biggl( \int_0^T
Q^{(n)}(r)\,\mathrm{d}r > M \biggr) = \pr \biggl(\int_0^{Tn}
Q(r) \,\mathrm{d}r > M n^2 \biggr).
\end{equation}
\end{lem}

\begin{pf}
Observe that the left-hand side of (\ref{eqlem2}) equals
\begin{eqnarray*}
&&
\pr \Biggl( \frac{1}{n} \int_0^T
\sup_{s \le r} \Biggl( \sum_{i = 1}^n
B^{(i)}(r) - B^{(i)}(s) - cn(r - s) \Biggr)\,\mathrm{d}r > M
\Biggr)
\\
&&\quad = \pr \biggl( \frac{1}{n} \int_0^T
\sup_{s \le r} \bigl( B(rn) - B(sn) - cn(r - s) \bigr)\,\mathrm{d}r > M
\biggr)
\\
&&\quad = \pr \biggl( \frac{1}{n} \int_0^T
\sup_{s \le rn} \bigl( B(rn) - B(s) - crn + cs \bigr)\,\mathrm{d}r > M \biggr).
\end{eqnarray*}
Using the substitution $v: = rn$, we obtain that
\[
\pr \biggl( \int_0^T Q^{(n)}(r)
\,\mathrm{d}r > M \biggr) = \pr \biggl( \int_0^{Tn}
\sup_{s \le v} \bigl(B(v) - B(s) - cv + cs\bigr) \,\mathrm{d}v > M
n^2 \biggr).
\]
This completes the proof.
\end{pf}

In our analysis, we use the following notation:
\[
\psi(M,a,s) := \frac{(M+({1}/{2})cs^2-as)^2}{({2}/{3})s^3}+2ac.
\]
The proof of Theorem \ref{thinterm} is based on the following lemmas.

%le3 #&#
\begin{lem}\label{comp}
For each $M,T >0$ it holds that
\[
\inf_{a\ge0} \inf_{s\in(0,T]} \psi(M,a,s) = \varphi(T,M).
\]
The optimizing $(a,s)$ equals
\[
(a\s,s\s)=\cases{ %
 (0,\sqrt{6M/c}),&\quad  $\mbox{{if} $\sqrt{6M/c}<T$};$
\vspace*{2pt}\cr
(M/T -cT/6,T), & \quad$\mbox{{else.}}$}
\]
\end{lem}

\begin{pf} Straightforward computation.
\end{pf}

Define
\[
p_n(T,M,a) := \mathbb P \biggl(\int
_0^T Q^{(n)}(r)\,\mathrm{d}r \ge M
\Big\vert Q^{(n)}(0)= a \biggr).
\]

%le4 #&#
\begin{lem}\label{bp}
For each $T, M ,a > 0$
\[
% \alpha(t,M,k,\epsilon)
% :=
\limsup_{n\to\infty}\frac{1}{n}\log
p_n(T,M,a) \le -\inf_{s\in[0,T]} \frac{(M+({1}/{2})cs^2 - a s)^2}{({2}/{3})s^3}.
\]
\end{lem}

\begin{pf}
The proof is based on the Schilder's sample-path large-deviations principle
\cite{Dembo98,Mandjes07}.
Define the path space
\[
\Omega:= \biggl\{f\dvtx \mathbb R \to\mathbb R, \mbox{continuous}, f(0) = 0,
\lim_{t \to\infty} \frac{f(t)}{1 + |t|} = \lim_{t \to-\infty}\frac
{f(t)}{1 + |t|} = 0
\biggr\},
\]
equipped with the norm
\[
\Vert f\Vert_\Omega:= \sup_{t \in\mathbb R} \frac{f(t)}{1 + |t|}.
\]
For a given function $f$, we have that the corresponding workload is
given through
$ q[f](t)
:=
\sup_{s\le t} (f(t)-f(s) - c(t-s)).
$
In addition,
\[
{\mathscr S} := \biggl\{f\in{\Omega}\dvt q[f](0)= a, \int_0^T
q[f](r) \,\mathrm{d}r \ge M \biggr\}.
\]
The set ${\mathscr S}$ is closed; the proof of this property can be
found in the \hyperref[app]{Appendix}.
Hence, due to Schilder's theorem, we have that
%
%e7 #&#
\begin{equation}
\label{ineqS} \limsup_{n\to\infty}\frac{1}{n}\log p_n(T,M,a)
\le -\inf_{f\in{\mathscr S}}{\mathbb I}(f),
\end{equation}
with
\[
{\mathbb I}(f) := \cases{ %
\displaystyle\frac{1}{2}\int
_\mathbb{R} \bigl(f'(r)\bigr)^2
\,\mathrm{d}r, & \quad $f \in\mathscr{A},$
\vspace*{2pt}\cr
\infty, &\quad  $\mbox{otherwise,}$}
\]
where
$\mathscr{A}
$
denotes the space of absolutely continuous functions with a square
integrable derivative.

Now we show that
%
%e8 #&#
\begin{equation}
\label{S=T} -\inf_{f \in\mathscr{S}}\mathbb{I}(f) = -\inf_{f \in\mathscr
{T}}
\mathbb{I}(f),
\end{equation}
where
\[
{\mathscr T} := \biggl\{ f\in\Omega\dvt \exists s \in[0,T]\dvt \int
_0^s f(r)\,\mathrm{d}r\ge M+\frac{1}{2}cs^2
- a s \biggr\}.
\]

To this end, first observe that $\mathscr{T} \subseteq\mathscr{S}$, so that
$-\inf_{f \in\mathscr{T}} \mathbb{I}(f) \le -\inf_{f \in\mathscr{S}}
\mathbb{I}(f)$; we are therefore
left with proving the opposite inequality.
Now fix for the moment a path $f$. Bearing in mind $f$ is an absolutely
continuous function, the
following procedure yields a path $\bar{f}\in{\mathscr T}$ with
${\mathbb I}(f)= {\mathbb I}(\bar f)$. First, we
let
\[
m[f] := \int_0^T 1_{\{q[f](u) > 0\}}(u)
\,\mathrm{d}u,
\]
denote the amount of `nonidle time' corresponding to the path $f$ in $[0,T]$.
Then define
\[
i[f](r) := \inf \biggl\{ s \in[0,T]\dvt \int_0^s
1_{\{q[f](u) > 0\}}(u) \,\mathrm{d}u > r \biggr\}
\]
for $r \in[0,m[f]]$,
and
\[
j[f](r):= \inf \biggl\{s \in[0,T]\dvt \int_0^s
1_{\{q[f](u) = 0\}}(u) \,\mathrm{d}u > r \biggr\}
\]
for
$r \in[0,T-m[f]]$.

Now we construct the path $\bar f$ by shifting all the idle periods of
$q[f]$ to the end of the interval $[0,T]$. That is,
for $r \in[0, m[f]]$, let
\[
\bar f(r) := q[f]\bigl(i[f](r)\bigr)+cr-a,
\]
and for $r \in[m[f], T]$, let
\[
\bar f(r) := q[f]\bigl(i[f]\bigl(m[f]\bigr)\bigr)+cr +f\bigl(j[f]\bigl(r-m[f]
\bigr)\bigr)-c j[f]\bigl(r-m[f]\bigr).
\]
We also set
\[
\bar{f}(r) := 0  \qquad \mbox{for }   r < 0\quad    \mbox{and} \quad   \bar{f}(r) :=
\bar{f}(T)\qquad   \mbox{for }   r >T.
\]
It is clear that ${\mathbb I}(f)= {\mathbb I}(\bar f)$ (because we just
permuted subintervals of $[0,T]$, which does not affect the rate
function), while the constructed path
$
\bar f
$ is now in ${\mathscr T}.$ Conclude that
$
- \inf_{f\in\mathscr{S}} \mathbb{I}(f)
\le
- \inf_{f\in\mathscr{T}} \mathbb{I}(f),
$ as desired.

We are therefore left with computing $- \inf_{f\in\mathscr{T}} \mathbb{I}(f)$.
Let $\epsilon> 0$.
Clearly, $
{\mathscr T}
\subseteq
\bigcup_{s \in[0,T]} {\mathscr T^s}
$, with
\[
\mathscr T^s := \biggl\{ f\in\Omega\dvt \int_0^s
f(r)\,\mathrm{d}r > M+\frac{1}{2}cs^2 - a s - \epsilon \biggr\}.
\]
This implies that
%
%e9 #&#
\begin{equation}
\label{inf} - \inf_{f \in\mathscr T} {\mathbb I}(f) \le - \inf_{s \in[0,t]}
\inf_{f \in\mathscr T^s} {\mathbb I}(f).
\end{equation}
Observe that set $\mathscr T^s$ is open, and combine this with
Schilder's theorem and \refs{intN}:
\begin{eqnarray*}
- \inf_{f \in\mathscr T^s} {\mathbb I}(f) %\nonumber
%&\le&
% - \inf_{f \in int \mathscr T^s} {\mathbb I}(f) \\
&\le& \lim_{n \to\infty} \frac{1}{n} \log \pr \Biggl(
\frac{1}{n} \int_0^s \sum
_{i = 1}^n B^{(i)}(r)\,\mathrm{d}r > M +
\frac{1}{2}c s^2 - as - \epsilon \Biggr)
\\
%&=&
% \lim_{n \to\infty} \frac{1}{n} \log
% \pr\left( \frac{1}{n} \sum_{i = 1}^n \mathcal{N}_i
% > M + \frac{1}{2}c s^2 - as \right) \\
&=&
\lim_{n \to\infty} \frac{1}{n} \log \pr \biggl( \mathscr{N} > \sqrt{
\frac{3n}{s^3}} \biggl(M + \frac{1}{2}c s^2 - as - \epsilon
\biggr) \biggr).
\end{eqnarray*}

%where ${\mathcal{N}_i}$ are independent normal random variables with
%mean $0$ and variance $\frac{1}{3} s^3$.
Using that
$\pr(\mathscr{N} > x) \le({\sqrt{2\uppi}x} )^{-1}\exp(-x^2/2)$,
we obtain
%
%e10 #&#
\begin{equation}
\label{ineqT} - \inf_{f\in{\mathscr T^s}}{\mathbb I}(f) \le - \frac{(M+({1}/{2})cs^2 - a s - \epsilon)^2}{({2}/{3})s^3}.
\end{equation}
%
%to this end realize that ${\mathbb V}{\rm ar} \int_0^s B(r){\rm d}r =
Thus the claim follows from combining (\ref{ineqS}), \refs{S=T}, and
(\ref{inf}) with
(\ref{ineqT}).
\end{pf}

\begin{pf*}{Proof of Theorem \ref{thinterm}}
Due to Lemma \ref{lmscaling} it suffices to find the logarithmic
asymptotics
\[
\lim_{n\toi}\frac{1}{n}\log\pr \biggl(\int_0^T
Q^{(n)}(r)\,\mathrm{d}r \ge M \biggr).
\]
We establish the
upper and lower bound separately.

\textit{Upper bound}:
Recall ${\mathbb P}(Q^{(n)}(0)\ge a) = \mathrm{e}^{-2nca}$ by virtue of (\ref
{distrQ}). For any $\epsilon>0$
and an arbitrary integer $N$,
\begin{eqnarray*}
 \pr \biggl(\int_0^T
Q^{(n)}(r)\,\mathrm{d}r \ge M \biggr) &=& \int_0^\infty2nc
\mathrm{e}^{-2ncv\epsilon} p_n(T,M,v\epsilon)\,\mathrm{d}v
\\
&\le& \sum_{k=0}^\infty 2nc
\mathrm{e}^{-2nck\epsilon}p_n\bigl(T,M,(k + 1)\epsilon\bigr)
\\
&\le& \sum_{k=0}^{N-1} 2nc \mathrm{e}^{-2nck\epsilon}
p_n\bigl(T,M,(k + 1)\epsilon\bigr)+ 2nc\cdot\frac
{\mathrm{e}^{-2ncN\epsilon}}{1-\mathrm{e}^{-2nc\epsilon}}.
\end{eqnarray*}
As a consequence, \cite{Dembo98}, Lemma 1.2.15, leads to
\begin{eqnarray*}
&& \limsup_{n\to\infty}\frac{1}{n}\log\pr \biggl(\int
_0^T Q^{(n)}(r)\,\mathrm{d}r \ge M
\biggr)
\\
&&\quad \le \max \biggl\{ \max_{k = 0,\ldots,N-1} \biggl( \lim_{n \to\infty}
\frac{1}{n} \log p_n\bigl(T,M,(k+1)\epsilon\bigr) - 2ck\epsilon
\biggr), -2cN\epsilon \biggr\}.
\end{eqnarray*}
Due to Lemma \ref{bp}, we can further bound this by
\begin{eqnarray*}
&&\max \Bigl\{ \max_{k = 0,\ldots,N-1} \Bigl(-\inf_{s\in[0,T]}\psi
\bigl(M,(k+1)\epsilon ,s\bigr) \Bigr)+2c\epsilon,-2cN\epsilon \Bigr\}
\\
&&\quad\le-\min \Bigl\{\inf_{a\ge0} \inf_{s\in[0,T]}\psi(M,a,s)-2c
\epsilon,2cN\epsilon \Bigr\}.
\end{eqnarray*}

% \alpha(t,M)
% \le
% -\min\left\{\min_{k\in\{0,\ldots,N\}}\left(\inf_{s\in[0,t]}
% -\psi(M,(k+1)\epsilon,s)\right) + 2ck\epsilon, 2c\epsilon(N+1)\right
Now Lemma \ref{comp} yields
\[
\limsup_{n\to\infty}\frac{1}{n}\log\pr \biggl(\int_0^T
Q^{(n)}(r)\,\mathrm{d}r \ge M \biggr) \le -\min \bigl\{\varphi(T,M)-2c
\epsilon,2c\epsilon N \bigr\}.
\]
We establish the upper bound by subsequently letting $N\uparrow
\infty$ and $\varepsilon\downarrow0.$

\textit{Lower bound}: Let $\epsilon>0$. Due to the Skorokhod
representation, we have, with $L^{(n)}(\cdot)$ defined in the
obvious way,
\[
Q^{(n)}(t) = Q^{(n)} (0) + \overline{B^{(n)}}(t) - ct
+ L^{(n)}(t).
\]
Observe that for each $a\ge0$ and $s\in[0,T]$,
\begin{eqnarray*}
\pr \biggl(\int_0^T
Q^{(n)}(r)\,\mathrm{d}r \ge M \biggr)
&\ge& \pr \biggl( \int
_0^s Q^{(n)}(r) \,\mathrm{d}r > M
\biggr)
\\
&\ge& \pr \biggl(\int_0^s
\bigl(Q^{(n)}(0) + \overline{B^{(n)}}(r) - cr\bigr) \,\mathrm{d}r
> M \biggr)
\\
&\ge& \int_a^{a+\varepsilon} 2nc\exp(-2ncv) \pr
\biggl( \int_0^s \overline{B^{(n)}}(r)
\,\mathrm{d}r > \frac{1}{2} c s^2 + M - vs \biggr) \,\mathrm{d}v
\\
&\ge& 2nc \exp\bigl(-2nc(a + \varepsilon)\bigr) \pr \biggl(\int
_0^s \overline{B^{(n)}}(r) \,\mathrm{d}r
> \frac{1}{2} c s^2 + M - as \biggr)
\\
&\ge& 2nc \exp\bigl(-2cn(a + \varepsilon)\bigr) \pr \biggl(
\mathscr{N} > \sqrt\frac{3n}{s^3} \biggl(\frac{1}{2} c s^2 +
M - as \biggr) \biggr).
\end{eqnarray*}
Now applying that for $x > 0$,
\[
\pr(\mathscr{N} > x) \ge\frac{x^2-1}{\sqrt{2\uppi}x^3} \exp\bigl(- x^2/2\bigr),
\]
see, for example, Section 2 in \cite{Adler90},
we obtain that for all $a\ge0$ and $s\in[0,T]$,
\[
\liminf_{n\to\infty}\frac{1}{n}\log\pr \biggl(\int_0^T
Q^{(n)}(r)\,\mathrm{d}r \ge M \biggr) \ge - \frac{ (({1}/{2}) c s^2 + M - as )^2}{({2}/{3})s^3} - 2c(a +
\varepsilon).
\]
In order to complete the proof it suffices to let $\epsilon
\downarrow0$ and to maximize over $a \ge0$ and $s \in[0,T]$.
\end{pf*}

%re1 #&#
\begin{rem}
Interestingly, the most likely path $f\s$ of
$\overline{B^{(n)}}(\cdot)$ can be explicitly computed, revealing
two separate scenarios.
\begin{itemize}
\item[--] Suppose $s\s=\sqrt{6M/c}<T$. Then the queue (most likely)
starts empty at time 0, is positive for a while, drops to 0 at time $s\s
$, and remains empty. The corresponding path $f\s$ of $\overline
{B^{(n)}}(\cdot)$ is, for $r\in[0,s\s]$,
\[
f\s(r) = 2cr -\frac{cr^2}{6}\sqrt{\frac{6c}{M}},
\]
and $f\s(r)=f\s(s\s)$ for $r\in(s\s,T].$
\item[--] Suppose $s\s=T<\sqrt{6M/c}$. Then the queue is symmetric in
the interval $[0,T]$, and has the value $a\s$ at times 0 and $T$.
The corresponding path $f\s$ of $\overline{B^{(n)}}(\cdot)$ is, for
$r\in[0,T]$,
\[
f\s(r) = 2cr-\frac{c}{T}r^2.
\]
\end{itemize}
It can easily be verified that indeed
\[
\frac{1}{2}\int_0^T \bigl((f
\s)'(r)\bigr)^2\,\mathrm{d}r + 2a\s c = \varphi(T,M)
\]
as expected.
\end{rem}

%s5 #&#
\section{Long timescale}\label{BMlong}
In this section, we consider the case that $T(u)$ is between $\sqrt{u}$
and $u$. It turns out that we find the same logarithmic asymptotics as
in the case that $T(u)=T\sqrt{u}$ for large $T$ (i.e., $T$ larger
than $\sqrt{6M/c}$).
In the proof, we first introduce some sort of `surrogate busy periods'
(recall that `traditional' busy periods do not exist for reflected
Brownian motion). Then we show that the event of interest occurs
essentially due to a single busy period being `big' (in terms of the
area swept under the workload graph); this is due to the fact that the
contribution of a single busy period has a subexponential distribution
(viz. roughly a Weibull distribution with shape parameter~$\frac{1}{2}$).

Defining
\[
\varphi(M) := \frac{2}{3} \sqrt{6} c \sqrt{cM}, \qquad    \tilde\psi(M,
\delta,s) := \frac{(M+({1}/{2})cs^2 - \delta s)^2}{({2}/{3})s^3},
\]
we are in a position to state the main result of the section.

%th3 #&#
\begin{thmm} \label{thlong}
Let $\sqrt{u} = \mathrm{o}(T(u)) $ and $T(u) = \mathrm{o}(u)$. Then,
\[
\lim_{u\to\infty}\frac{1}{\sqrt{u}}\log \pr \biggl(\int_0^{T(u)}
Q(r)\,\mathrm{d}r > Mu \biggr) = - \varphi(M).
\]
\end{thmm}

In order to prove Theorem \ref{thlong}, we need to introduce some
notation. Let
\[
\tau_0:= \inf\bigl\{t > 0\dvt Q(0) + B(t) - ct = 0\bigr\}, \qquad    \tau(x):= \inf\bigl\{t > 0\dvt x + B(t) - ct = 0\bigr\}.
\]

Besides, for {given $\delta>0$ and} $i = 1,2,\ldots,$ let
\[
\sigma_i := \inf\bigl\{t > \tau_{i-1}\dvt Q(t) \ge2 \delta
\bigr\},   \qquad  \tau_i := \inf\bigl\{t > \sigma_i\dvt
Q(t) \le\delta\bigr\}
\]
and
\[
H_0 := \int_0^{\tau_0} Q(r)
\,\mathrm{d}r,  \qquad    H_i := \int_{\sigma_i}^{\tau_i}
Q(r) \,\mathrm{d}r.
\]
Observe that $\{H_i\}_{i\in{\mathbb N}}$ constitutes a sequence of
i.i.d. random variables, that is in addition independent of $H_0$;
likewise, the $\xi_i:=\tau_i-\sigma_i$ are i.i.d. random variables.
Moreover, for each $i = 1,2,\ldots$ we have
\[
H_i \stackrel{\,\mathrm{d}} {=} \int_0^{\tau(\delta)}
\bigl(\delta+ B(r) - cr\bigr) \,\mathrm{d}r.
\]

The following lemmas play crucial role in the proof of Theorem \ref{thlong}.

%le5 #&#
\begin{lem} \label{lemP1P2}
For each $M > 0$ it holds that
\[
\lim_{\delta\downarrow0} \inf_{s \ge0} \tilde\psi(M,\delta,s) = \varphi(M).
\]
\end{lem}
\begin{pf} This proof is a straightforward computation.
Note that
\[
s^\star(\delta) = \frac{-\delta+ \sqrt{\delta^2 + 6Mc}}{c}
\]
is the minimizer in $\inf_{s \ge0} \tilde\psi(M,\delta,s)$. Consequently,
\[
\lim_{\delta\downarrow0} \inf_{s \ge0} \tilde\psi(M,\delta,s) =
\lim_{\delta\downarrow0} \frac{(M+({1}/{2})c(s^\star(\delta))^2 - \delta s^\star(\delta
))^2}{({2}/{3})(s^\star(\delta))^3} = \varphi(M).
\]
This completes the proof.
\end{pf}

%le6 #&#
\begin{lem} \label{lemlong}
For each $M > 0$ and $i = 0, 1,\ldots,$ we have
\[
\limsup_{u \toi} \frac{1}{\sqrt{u}} \log \pr (H_i > Mu ) \le
-\varphi(M).
\]
\end{lem}

\begin{pf}
We start with the analysis of $H_i$, for $i = 1,2,\ldots.$
Observe that
\[
\pr (H_i > Mu ) = \pr \biggl(\exists s \ge0 \dvt \frac{1}{u}
\int_0^s \bigl(\delta+ B(r) - cr\bigr)
\,\mathrm{d}r > M, \forall r \in(0,s) \dvt \delta+ B(r) - cr > 0 \biggr),
\]
which is majorized by
%
%e11 #&#
\begin{equation}
\label{eqlongH1} \pr \biggl(\exists s \ge0 \dvt \frac{1}{u} \int
_0^s \bigl(\delta+ B(r) - cr\bigr) \,\mathrm{d}r
> M \biggr).
\end{equation}
Substituting $r = \sqrt{u}v$ we obtain that, for $u$ sufficiently
large, \refs{eqlongH1} equals
%
%e12 #&#
\begin{eqnarray}\label{eqlongH2}
&& \pr \biggl(\exists s \ge0 \dvt \int_0^{{s}/{\sqrt{u}}} \biggl(\frac{\delta}{\sqrt{u}} + \frac{1}{\sqrt{u}}B(\sqrt{u}v) - c v
\biggr) \,\mathrm{d}v > M \biggr)
\nonumber\\
&&\quad= \pr \biggl(\exists s \ge0 \dvt \int_0^s
\biggl(\frac{\delta}{\sqrt{u}} + u^{-1/4}B(v) - c v \biggr) \,\mathrm{d}v > M
\biggr)
\\
&&\quad=  \pr \biggl(\sup_{s \ge0} { \frac{ \int_0^s B(v) \,\mathrm{d}v}{ M + ({1}/{2}) cs^2 - {\delta s}/{\sqrt{u}}}} >
u^{1/4} \biggr).\nonumber
\end{eqnarray}
Now, observe that
$Y(s) := \int_0^s B(v)\,\mathrm{d}v / (M + \frac{1}{2} cs^2 -{\delta
s}/{\sqrt{u}} )$ has bounded trajectories a.s. Hence, the Borell
inequality (see, e.g., \cite{Adler90}, Theorem 2.1) leads to the
following upper bound of \refs{eqlongH2}:
\[
2\exp \biggl(- \inf_{s \ge0}\frac{(M+({1}/{2})cs^2 - ({\delta}/{\sqrt{u}}) s)^2}{({2}/{3})s^3} \Bigl(u^{{1}/{4} }- \E\sup_{s \ge0} Y(s) \Bigr)^2 \biggr),
\]
where
$
\E\sup_{s \ge0} Y(s)
$
is bounded (by `Borell').
Combining the above with Lemma \ref{lemP1P2}, we obtain that
%
%e13 #&#
\begin{equation}
\label{eqHiupper} \limsup_{u \toi} \frac{1}{\sqrt{u}} \log\pr
(H_i > Mu ) \le - \varphi(M).
\end{equation}
%
%In order to obtain (i) it suffices to observe that $\tilde\psi(

In order to prove the claim for $H_0$ observe that
\[
\pr (H_0 > Mu ) = \int_0^\infty2ca
\exp(-2ca) \pr \biggl(\int_0^{\tau(a)} \bigl(a + B(r)
- cr\bigr) \,\mathrm{d}r > Mu \biggr) \,\mathrm{d} a.
\]
Thus, by \refs{eqHiupper}, it suffices to proceed along the lines
of the proof of the upper bound of Theorem~\ref{thinterm}.
\end{pf}

\begin{pf*}{Proof of Theorem \ref{thlong}}
We establish upper and lower bound separately.

\textit{Lower bound}:
The lower bound follows straightforwardly from Theorem \ref{thinterm}
combined with the fact that
for sufficiently large $u$ we have (recalling that $\sqrt{u}= \mathrm{o}(T(u))$)
\[
\pr \biggl(\int_0^{T(u)} Q(r)\,\mathrm{d}r > Mu
\biggr) \ge \pr \biggl(\int_0^{\sqrt{({6M}/{c}) u}} Q(r)
\,\mathrm{d}r > Mu \biggr).
\]

\textit{Upper bound}:
Let $\delta>0$ and denote
$
N(u):= \inf\{i \dvt \tau_i \ge T(u)\} $,
$
K := 2/\E\xi_i.
$
Observe that
\begin{eqnarray*}
 \pr \biggl(\int_0^{T(u)} Q(r)
\,\mathrm{d}r > u \biggr) &\le& \pr \Biggl(2\delta T(u) + \sum
_{i = 0}^{N(u)} H_i > u \Biggr)
\\
&\le& \pr \Biggl(2\delta T(u) + \sum_{i = 0}^{N(u)}
H_i > u, N(u) \le K T(u) \Biggr) + \pr \bigl(N(u) > K T(u) \bigr)
\\
&\le&\bar{P}_1(u) + \bar{P}_2(u),
\end{eqnarray*}
with
\[
\bar{P}_1(u):= \pr \Biggl(\sum_{i = 0}^{\lceil K T(u)\rceil}
H_i > u - 2\delta T(u) \Biggr) ,\qquad  \bar{P}_2(u):= \pr
\bigl(N (u) > \bigl\lfloor KT(u)\bigr\rfloor\bigr).
\]

We first analyze $\bar{P}_1(u)$.
The idea is to reduce the problem of finding the upper bound of $\bar{P}_1(u)$
to the setting of \cite{Denisov08}, Theorem 8.3.
To this end, pick $\varepsilon> 0$.
Due to Lemma \ref{lemlong} there exists a sequence $\{\tilde{H}_i\}_{i=0,1,\ldots}$ of i.i.d. random variables such that for each $x > 0$
and $\delta$ sufficiently small,
%
%e14 #&#
\begin{equation}
\label{eqHtilde} \pr(H_i > x) \le\pr(\tilde{H}_i > x)
\end{equation}
and
%
%e15 #&#
\begin{equation}
\label{eqHias} \pr(\tilde{H}_i > x) = p(x) \exp\bigl(-\bigl(
\varphi(M) - \varepsilon\bigr) \sqrt{x}\bigr),
\end{equation}
where $p(\cdot)$ is some $O$-regularly varying function, that is, $p(x)$
is a
measurable function, such that, for each $\lambda\ge1$
\[
0< \liminf_{x \toi} \frac{p(\lambda x)}{p(x)} \le\limsup_{x \toi}
\frac
{p(\lambda x)}{p(x)} < \infty
\]
(see, e.g., \cite{Bingham87}, Chapter 2, or the Appendix of \cite
{Denisov08}).

It is standard that, due to \refs{eqHtilde}, for each $x > 0$,
%
%e16 #&#
\begin{equation}
\pr \Biggl(\sum_{i=0}^{\lceil K T(u)\rceil}
H_i > x \Biggr) \le \label{eqHisum} \pr \Biggl(\sum
_{i=0}^{\lceil K T(u)\rceil} \tilde{H}_i > x \Biggr).
\end{equation}
Now, applying \cite{Denisov08}, Theorem 8.3 and recalling that $KT(u) =
\mathrm{o}(u)$, we {have}, as $u\toi$,
%
%e17 #&#
\begin{equation}
\pr \Biggl(\sum_{i = 0}^{\lceil K T(u)\rceil}
\tilde{H}_i > u - 2\delta T(u) \Biggr) = \label{eqHisumas} \bigl
\lceil KT(u)\bigr\rceil\cdot\pr\bigl(\tilde{H}_0 > u - 2\delta T(u)
\bigr) \bigl(1 + \mathrm{o}(1)\bigr).
\end{equation}
Combining \refs{eqHisum} and \refs{eqHisumas} with \refs{eqHias},
we obtain that, for each $\varepsilon> 0$,
\[
\limsup_{u \toi} \frac{1}{\sqrt{u}} \log\bar{P}_1(u) \le -
\varphi(M) + \varepsilon;
\]
%
%Now it suffices to pass with $\varepsilon\to0$.
letting $\varepsilon\downarrow0$, we conclude that we can replace
the right-hand side in the
previous display by $-\varphi(M).$

We now focus on $\bar{P}_2(u)$.
Observe that
\[
\bar{P}_2(u) \le \pr\bigl(S_{\lfloor KT(u) \rfloor} < T(u)\bigr)\qquad
\mbox{where } S_{\lfloor K T(u) \rfloor} := \tau_0 + \sum
_{i = 1}^{\lfloor KT(u)
\rfloor} \xi_i.
\]
Moreover, note that $\xi_i, i = 1,2,\ldots$ are i.i.d. with
\[
\frac{\mathrm{d}}{\mathrm{d}t}\pr(\xi_1 \le t) = \frac{\delta}{\sqrt{2\uppi t^3}} \exp
\bigl(-(\delta- ct)^2/2t\bigr)
\]
{for $t>0$; see, for example, \cite{Rogers00}, Section 2.9.}
Hence,
%there exists $\theta_0 < 0$ such that
a Chernoff bound argument yields, recalling that $K > 1/\E\xi_i$,
\[
\limsup_{u \toi} \frac{1}{T(u)} \log\pr\bigl(S_{\lfloor KT(u)\rfloor} < T(u)
\bigr) \le - K\cdot\sup_{\theta< 0} \biggl[\theta\frac{1}{K} - \log\E
\exp (\theta\xi_1) \biggr] < 0.
\]
We have found that $\bar{P}_1(u)$ is smaller than a function of the
order $\exp(-\beta_1\sqrt{u})$, while
$\bar{P}_2(u)$ is smaller than a function of the order $\exp(-\beta_2
T(u))$, for some $\beta_1,\beta_2>0$. Now recalling that $\sqrt {u}=\mathrm{o}(T(u))$,
it follows that the upper bound on $\bar{P}_1(u)$ is
smaller than the upper bound on $\bar{P}_2(u)$.
As a result,
\[
\limsup_{u \toi} \frac{1}{T(u)} \log\pr \biggl(\int
_0^{T(u)} Q(r) \,\mathrm{d}r > Mu \biggr) \le
\limsup_{u\toi}\frac{1}{\sqrt{u}}\log\bigl(\bar{P}_1(u)+
\bar{P}_2(u)\bigr) = -\varphi(M).
\]
This completes the proof.
\end{pf*}

%s6 #&#
\section{Residual busy period} \label{BMrunBP}
\newcommand{\AI}{\mathrm{Ai}}

In this section, we analyze the integral of the stationary workload
for regulated Brownian motion over the \textit{residual busy period}.
It turns out to be possible to explicitly compute its Laplace transform,
in terms of the so-called Airy function.
As a by-product, the corresponding mean value is calculated.
%, which also enables us to find an expression for all moments.
%In particular, the corresponding mean and variance are calculated.

Recall that
\[
\tau_0:=\inf\bigl\{t \ge0 \dvt Q(t) = 0 \bigr\}, \qquad     \tau(x)
:= \inf\bigl\{t \ge0 \dvt x + B(t) - ct = 0 \bigr\};
\]
we also define the integral of the workload until the end of the busy
period, conditional on
the workload being $x$ at time 0:
\[
J(x):=\int_0^{\tau(x)} \bigl(x + B(t) - ct\bigr)
\,\mathrm{d}t.
\]
By
\[
\AI(x) := \frac{1}{\uppi} \int_0^\infty\cos
\biggl(\frac{1}{3} t^3 + xt \biggr) \,\mathrm{d}t
\]
we denote the \textit{Airy function} (see, e.g., \cite{Abramovitz72}, Chapter 10.4).

%th4 #&#
\begin{thmm}\label{thbusyperiod}
For each $\gamma\ge0$,
\[
\E\exp \biggl[- \gamma\int_0^{\tau_0} Q(t)
\,\mathrm{d}t \biggr] = \frac{2c}{\AI ( (2\gamma)^{-{2}/{3}} c^2  )} \int_0^\infty
\mathrm{e}^{-cx} \AI \bigl( (2\gamma)^{-{2}/{3}} c^2 + (2
\gamma)^{{1}/{3}} x \bigr)\,\mathrm{d}x.
\]
\end{thmm}

\begin{pf}
Observe that up to time {$\tau_0$} we have that $Q(t) =Q(0) + B(t) -
ct.$ Hence,
%
%e18 #&#
\begin{eqnarray}\label{LTll}
&&\E\exp \biggl[- \gamma\int_0^{\tau_0}
Q(t) \,\mathrm{d}t \biggr]
\nonumber\\
&&\quad= \int_0^\infty\pr \biggl( \exp \biggl[ -\gamma
\int_0^{\tau_0} Q(t) \,\mathrm{d}t \biggr] > u \biggr)
\,\mathrm{d}u
\nonumber
\\
&&\quad= \int_0^\infty\pr \biggl( \exp
\biggl[ -\gamma\int_0^{\tau_0} \bigl(Q(0) + B(t) - ct
\bigr) \,\mathrm{d}t \biggr] > u \biggr) \,\mathrm{d}u
\nonumber
\\[-8pt]
\\[-8pt]
\nonumber
&&\quad= \int_0^\infty\int
_0^\infty2c \exp(-2cx) \pr \bigl( \exp \bigl[ -\gamma
J(x) \bigr] > u \bigr) \,\mathrm{d}x \,\mathrm{d}u
\\
\nonumber
&&\quad= \int_0^\infty2c \exp(-2cx) \int
_0^\infty \pr \bigl( \exp \bigl[ -\gamma J(x)
\bigr] > u \bigr) \,\mathrm{d}u \,\mathrm{d}x
\\
&&\quad=  \int_0^\infty2c \exp(-2cx) \E
\bigl[ \exp \bigl[ -\gamma J(x) \bigr] \bigr] \,\mathrm{d}x.\nonumber
\end{eqnarray}

Following Borodin and Salminen \cite{Borodin02}, Chapter 2, equation (2.8.1), we
have that
%
%e19 #&#
\begin{equation}
\label{eqBS} \E \bigl[ \exp \bigl[ -\gamma J(x) \bigr] \bigr] = \exp(cx)
\frac{\AI (2^{{1}/{3}} \gamma^{-{2}/{3}} (({1}/{2}) c^2 +
\gamma x)  )} {
\AI ( (2\gamma)^{-{2}/{3}} c^2  )},
\end{equation}
which combined with \refs{LTll} completes the proof.
\end{pf}

In the following proposition, we compute the mean value of the integral
over the residual busy
period, given the workload at time 0 equals $x$.

%pr1 #&#
\begin{prop}\label{ptransient}
The mean area until the end of the transient busy period, is
\[
\E J(x)= \E \biggl[ \int_0^{\tau(x)} \bigl(x + B(t)
- ct\bigr) \,\mathrm{d}t \biggr] = \frac{x^2}{2c} + \frac{x}{2c^2}.
\]
\end{prop}

\begin{pf}
Due to the fact that
%
%e20 #&#
\begin{equation}
\label{Airyinfinity}
\AI(u) = \frac{1}{2\sqrt{\uppi} u^{{1}/{4}}} \exp \biggl(-\frac{2}{3}
u^{{3}/ {2}} \biggr) \biggl(1 -\frac{5}{48} u^{-3/2}+\mathrm{o}
\bigl(u^{-3/2}\bigr) \biggr)
\end{equation}
as $u\to\infty$, combined with \refs{eqBS}, we have that
\begin{eqnarray*}
\E \bigl[ - \gamma J(x) \bigr] &=& \exp(cx) \biggl( \frac{1}{1 + ({2 x}/{c^2}) \gamma}
\biggr)^{{1}/{4}} \exp \biggl[\frac{c^3}{3\gamma} \biggl(1 - \biggl(1 +
\frac{2\gamma x}{c^2} \biggr)^{{3}/ {2}} \biggr) \biggr]\bigl(1+\mathrm{o}(\gamma)
\bigr)
\\
% &=&
% \exp(cx)
% \left( \frac{1}{1 + \frac{2 x}{c^2} \gamma} \right)^\frac{1}{4}
% \exp\left(\frac{2}{3} (2\gamma)^{-1} c^3 -
% \frac{2}{3}\left((2\gamma)^{-\frac{2}{3}} c^2 + t
% (2\gamma)^\frac{1}{3} x\right)^{3/2} \right) \\
&=& 1 -
\biggl(\frac{x^2}{2c} + \frac{x}{2c^2} \biggr) \gamma+ \mathrm{o}(\gamma)
\end{eqnarray*}
as $\gamma\to0$.
This completes the proof.\vadjust{\goodbreak}
\end{pf}

Combining Proposition \ref{ptransient} with
\refs{LTll} (and using the dominated convergence theorem)
immediately leads to the following corollary.
%
%co1 #&#
\begin{cor} \label{mean}
\[
\E \biggl[ \int_0^{\tau_0} Q(t) \,\mathrm{d}t \biggr]
= \frac{1}{2c^3}.
\]
\end{cor}

We note that, applying more precise expansions in \refs{Airyinfinity},
one can get the analogue of Proposition \ref{ptransient}
for higher moments of $J(x)$,
and (by applying \refs{LTll}) also formulas for corresponding moments of
$\int_0^{\tau_0} Q(t) \,\mathrm{d}t$. These computations are tedious
(although standard),
and are therefore left out.

%s7 #&#
\section{Discussion and outlook}
In this paper, we analyzed the probability that the area swept under
{the Brownian storage graph}
between $0$ and $T(u)$ exceeds $u$.
{We did so for various types of interval lengths $T(u)$, leading to
asymptotic results for three timescales ($u\toi$).
{A topic for future research could be to consider a wider class of
inputs} $\{X(t)\dvt t\in{\mathbb R}\}$,
for instance, Gaussian processes or L\'evy processes. In the former case,
there is the major complication that $Q(0)$ is \textit{not} independent
of $\{X(t)\dvt t>0\}$,
which is a property that we repeatedly used in this paper.
In the latter case, we have to make sure that all steps in which we use specific
properties of Brownian motion, carry over to the more general L\'evy case.
We do anticipate, though, that in case the L\'evy-input is light-tailed
the asymptotics are in the qualitative sense very similar to those
related to the Brownian case (i.e., the same three regimes apply).
Another related problem concerns the derivation of a central limit
theorem for
\[
\frac{1}{\sqrt{T}} \biggl(\int_0^T Q(t)
\,\mathrm{d}t - q{T} \biggr),
\]
with $q$ the mean stationary workload.

%s8 #&#
\begin{appendix}\label{app}
\section*{Appendix}
In this appendix, we prove that
\[
\mathscr S = \biggl\{f \in\Omega\dvt q[f](0) = a, \int_0^T
q[f](s) \,\mathrm{d}s \ge M \biggr\}
\]
is a closed set in the space $\Omega$. To this end, let $f_n \in
\mathscr S$ be a sequence of functions such that $\Vert f_n -
f\Vert_\Omega\to0$, as $n \to\infty$ for some function $f \in
\Omega$. We prove our claim by showing that $f \in\mathscr S$.

First, we show that for the limiting path $f$ it holds that
%
%e21 #&#
\renewcommand{\theequation}{\arabic{equation}}
\setcounter{equation}{19}
\begin{equation}
\label{q0} q[f](0) = - \inf_{s\le0} \bigl(f(s) - cs\bigr) = a.
\end{equation}
%
%Observe that for each $f \in\mathscr S$ there exists $s_f$ such that
%$\inf_{s\le0} f(s) = f(s_f)$.
First, observe that $g(s) - cs \to\infty$ as $s \to-\infty$, as an
immediate consequence of the fact that ${|g(s) - cs|}/({1 + |s|})
\to c$ for all $g \in\Omega$. Consequently, for any such $g$ there
is a point $s$ in which~$g$ takes its minimum in $[-\infty, 0]$.

Let $s_0$ be such that $ \inf_{s\le0} (f(s) - cs) = f(s_0) - c s_0
$. Then
\[
- a \le \lim_{n \to\infty} f_n(s_0) - c
s_0 = f(s_0) -c s_0.
\]

On the other hand, let $\{s_n\}$ be the sequence of points such that
$
\inf_{s\le0} (f_n(s) - cs) = f(s_n) - c s_n.
$ Observe that $\{s_n\}$ is bounded. If not, then, for each $k$ and
$\varepsilon> 0$, we would have
\[
\bigl\Vert f_k(s) - f(s)\bigr\Vert_\Omega \ge \sup_{ s \in\{s_n\}}
\frac{|f_k(s) - f(s)|}{1 + |s|} = \sup_{s \in\{s_n\}} \frac{|-a - f(s) - cs|}{1 + |s|} \ge c - \varepsilon.
\]
Conclude that there exists an $M >0$ such that $|s_n| < M$. For $n$
large enough
\[
\bigl|f(s_n) - c s_n - \bigl(f_n(s_n)
- c s_n\bigr)\bigr| = \bigl|f_n(s_n) -
f(s_n)\bigr| \le \bigl(1 + |s_n|\bigr) \varepsilon \le (1 + M)
\varepsilon,
\]
which implies
\[
f(s_0) - c s_0 \le f(s_n) -
cs_n \le f_n(s_n) - c s_n + (1 +
M) \varepsilon = - a + (1 + M) \varepsilon.
\]
To complete the proof of (\ref{q0}), it is enough to let
$\varepsilon\downarrow0$.

Now we prove that
\[
\int_0^T q[f](s)\,\mathrm{d}s \ge M.
\]

Observe that
\[
\int_0^T \bigl|q[f_n](s) - q[f](s)\bigr|
\,\mathrm{d}s \le I_1+I_2,
\]
where
\[
I_1:= \int_0^T
\bigl|f_n(s) - f(s)\bigr| \,\mathrm{d}s, \qquad I_2:= \int
_0^T \Bigl\llvert \inf_{r\le s}\bigl(
f(r) - cr\bigr) - \inf_{v\le s} \bigl(f_n(v) - cv\bigr)\Bigr
\rrvert \,\mathrm{d}s.
\]
Let us examine $I_1$ first. Due to the fact that
$\lim_{n \to\infty} \Vert f_n - f\Vert_\Omega= 0,$ we have for $n$ large
enough that
%
%e22 #&#
\begin{equation}
\label{sup} \varepsilon \ge \sup_{s\le T} \frac{|f_n(s) - f(s)|}{1 + |s|} \ge
\sup_{s\in[0,T]} \frac{|f_n(s) - f(s)|}{1 + s} \ge \sup_{s\in[0,T]} \frac{|f_n(s) - f(s)|}{1 + T}.
\end{equation}
This implies
\[
\int_0^T \bigl|f_n(s) - f(s)\bigr|
\,\mathrm{d}s < T(1+T) \varepsilon.
\]

Now consider $I_2$. Let $s_0$ be the minimizer in $\inf_{r \in
[0,s]} (f(r) - cr)$ and $s_n$ the minimizer in $\inf_{r \in[0,s]}
(f_n(r_n) - cr_n)$. Then (\ref{sup}) implies that for $n$ large
enough
\[
f_n(s_n) - c s_n - \bigl(f(s_0)
- c s_0\bigr) \le f_n(s_0) - c
s_0 - \bigl(f(s_0) - c s_0\bigr) \le (1 +
T)\varepsilon.
\]
On the other hand
\[
f(s_0) - c s_0 - \bigl(f_n(s_n)
- cs_n\bigr) \le f(s_n) - c s_n -
\bigl(f_n (s_n) - c s_n\bigr) \le (1 + T)
\varepsilon.
\]
It follows that $I_2\le T(1+T)\varepsilon.$ Now it is enough to let
$\varepsilon\downarrow0$; realizing that for each $n$ we have $
\int_0^T q[f_n](s) \,\mathrm{d}s \ge M,
$ the proof is completed.
\end{appendix}

\section*{Acknowledgements}
K. D\c{e}bicki and M. Mandjes thank the Isaac Newton Institute,
Cambridge, for hospitality.
Jose Blanchet (Columbia University, New York), Peter Glynn (Stanford
University),
Sean Meyn (Univ. of Illinois at Urbana-Champaign), and Florian Simatos
(CWI, Amsterdam)
are thanked for valuable comments and inspiring discussions.

M. Arendarczyk was supported by MNiSW Grant N N201 412239 (2010--2011),
and K.~D\c{e}bicki by MNiSW Grant N N201 394137 (2009--2011).

% imsref loaded by akundreckaite, 2013-01-11 12:12:31
%

%suskaldyti doi

\printhistory


\begin{thebibliography}{16}
% BibTex style file: bej.bst, 2012-09-27
% Default style options (sort=1,type=number).
% Used options (sort=1,type=number).

%b1 #&#
\bibitem{Abramovitz72}
%
\begin{bbook}[auto:STB|2013/01/04|07:53:31]
\bauthor{\bsnm{Abramovitz},~\bfnm{M.}\binits{M.}} \AND
\bauthor{\bsnm{Stegun},~\bfnm{I.}\binits{I.}}
(\byear{1972}).
\btitle{Handbook of Mathematical Functions with Formulas, Graphs, and
Mathematical Tables}.
\blocation{New York, NY}: \bpublisher{Wiley}.
\bptok{imsref}%
\end{bbook}
%
\endbibitem

%b2 #&#
\bibitem{Adler90}
%
\begin{bbook}[mr]
\bauthor{\bsnm{Adler},~\bfnm{Robert~J.}\binits{R.J.}}
(\byear{1990}).
\btitle{An Introduction to Continuity, Extrema, and Related Topics for General
{G}aussian Processes}.
\bseries{Institute of Mathematical Statistics Lecture Notes---Monograph
Series}
\bvolume{12}.
\blocation{Hayward, CA}: \bpublisher{IMS}.
\bid{mr={1088478}}
\bptok{imsref}%
\end{bbook}
%
\endbibitem

%b3 #&#
\bibitem{Asmussen07}
%
\begin{bbook}[mr]
\bauthor{\bsnm{Asmussen},~\bfnm{S{\o}ren}\binits{S.}} \AND
\bauthor{\bsnm{Glynn},~\bfnm{Peter~W.}\binits{P.W.}}
(\byear{2007}).
\btitle{Stochastic Simulation: Algorithms and Analysis}.
\bseries{Stochastic Modelling and Applied Probability}
\bvolume{57}.
\blocation{New York}: \bpublisher{Springer}.
\bid{mr={2331321}}
\bptok{imsref}%
\end{bbook}
%
\endbibitem

%b4 #&#
\bibitem{Bingham87}
%
\begin{bbook}[mr]
\bauthor{\bsnm{Bingham},~\bfnm{N.~H.}\binits{N.H.}},
\bauthor{\bsnm{Goldie},~\bfnm{C.~M.}\binits{C.M.}} \AND
\bauthor{\bsnm{Teugels},~\bfnm{J.~L.}\binits{J.L.}}
(\byear{1987}).
\btitle{Regular Variation}.
\bseries{Encyclopedia of Mathematics and Its Applications}
\bvolume{27}.
\blocation{Cambridge}: \bpublisher{Cambridge Univ. Press}.
\bid{mr={0898871}}
\bptok{imsref}%
\end{bbook}
%
\endbibitem

%b5 #&#
\bibitem{BlanchetGlynnMeyn11}
%
\begin{bmisc}[auto:STB|2013/01/04|07:53:31]
\bauthor{\bsnm{Blanchet},~\bfnm{J.}\binits{J.}},
\bauthor{\bsnm{Glynn},~\bfnm{P.}\binits{P.}} \AND
\bauthor{\bsnm{Meyn},~\bfnm{S.}\binits{S.}}
(\byear{2011}).
\bhowpublished{Large deviations for the empirical mean of an M$/${M}$/$1 queue.
Unpublished manuscript.}
\bptok{imsref}%
\end{bmisc}
%
\endbibitem

%b6 #&#
\bibitem{Borodin02}
%
\begin{bbook}[mr]
\bauthor{\bsnm{Borodin},~\bfnm{Andrei~N.}\binits{A.N.}} \AND
\bauthor{\bsnm{Salminen},~\bfnm{Paavo}\binits{P.}}
(\byear{2002}).
\btitle{Handbook of {B}rownian Motion---Facts and Formulae},
\bedition{2nd} ed.
\bseries{Probability and Its Applications}.
\blocation{Basel}: \bpublisher{Birkh\"auser}.
\bid{doi={10.1007/978-3-0348-8163-0}, mr={1912205}}
\bptok{imsref}%
\end{bbook}
%
\endbibitem

%b7 #&#
\bibitem{Borovkov03a}
%
\begin{barticle}[mr]
\bauthor{\bsnm{Borovkov},~\bfnm{A.~A.}\binits{A.A.}},
\bauthor{\bsnm{Boxma},~\bfnm{O.~J.}\binits{O.J.}} \AND
\bauthor{\bsnm{Palmowski},~\bfnm{Z.}\binits{Z.}}
(\byear{2003}).
\btitle{On the integral of the workload process of the single server queue}.
\bjournal{J. Appl. Probab.}
\bvolume{40}
\bpages{200--225}.
\bid{issn={0021-9002}, mr={1953775}}
\bptok{imsref}%
\end{barticle}
%
\endbibitem

%b8 #&#
\bibitem{Dembo98}
%
\begin{bbook}[mr]
\bauthor{\bsnm{Dembo},~\bfnm{Amir}\binits{A.}} \AND
\bauthor{\bsnm{Zeitouni},~\bfnm{Ofer}\binits{O.}}
(\byear{1998}).
\btitle{Large Deviations Techniques and Applications},
\bedition{2nd} ed.
\bseries{Applications of Mathematics (New York)}
\bvolume{38}.
\blocation{New York}: \bpublisher{Springer}.
\bid{mr={1619036}}
\bptok{imsref}%
\end{bbook}
%
\endbibitem

%b9 #&#
\bibitem{Denisov08}
%
\begin{barticle}[mr]
\bauthor{\bsnm{Denisov},~\bfnm{D.}\binits{D.}},
\bauthor{\bsnm{Dieker},~\bfnm{A.~B.}\binits{A.B.}} \AND
\bauthor{\bsnm{Shneer},~\bfnm{V.}\binits{V.}}
(\byear{2008}).
\btitle{Large deviations for random walks under subexponentiality: The big-jump
domain}.
\bjournal{Ann. Probab.}
\bvolume{36}
\bpages{1946--1991}.
\bid{doi={10.1214/07-AOP382}, issn={0091-1798}, mr={2440928}}
\bptok{imsref}%
\end{barticle}
%
\endbibitem

%b10 #&#
\bibitem{Duffy10}
%
\begin{barticle}[auto:STB|2013/01/04|07:53:31]
\bauthor{\bsnm{Duffy},~\bfnm{K.}\binits{K.}} \AND
\bauthor{\bsnm{Meyn},~\bfnm{S.}\binits{S.}}
(\byear{2010}).
\btitle{Most likely paths to error when estimating the mean of a reflected
random walk}.
\bjournal{Perf. Eval.}
\bvolume{67}
\bpages{1290--1303}.
\bptok{imsref}%
\end{barticle}
%
\endbibitem

%b11 #&#
\bibitem{Kyprianou06}
%
\begin{bbook}[mr]
\bauthor{\bsnm{Kyprianou},~\bfnm{Andreas~E.}\binits{A.E.}}
(\byear{2006}).
\btitle{Introductory Lectures on Fluctuations of {L}\'evy Processes with
Applications}.
\bseries{Universitext}.
\blocation{Berlin}: \bpublisher{Springer}.
\bid{mr={2250061}}
\bptok{imsref}%
\end{bbook}
%
\endbibitem

%b12 #&#
\bibitem{Mandjes07}
%
\begin{bbook}[mr]
\bauthor{\bsnm{Mandjes},~\bfnm{Michel}\binits{M.}}
(\byear{2007}).
\btitle{Large Deviations for {G}aussian Queues: Modelling Communication
Networks}.
\blocation{Chichester}: \bpublisher{Wiley}.
\bid{doi={10.1002/9780470515099}, mr={2329270}}
\bptok{imsref}%
\end{bbook}
%
\endbibitem

%b13 #&#
\bibitem{Meyn08}
%
\begin{bbook}[mr]
\bauthor{\bsnm{Meyn},~\bfnm{Sean}\binits{S.}}
(\byear{2008}).
\btitle{Control Techniques for Complex Networks}.
\blocation{Cambridge}: \bpublisher{Cambridge Univ. Press}.
\bid{mr={2372453}}
\bptok{imsref}%
\end{bbook}
%
\endbibitem

%b14 #&#
\bibitem{Meyn06}
%
\begin{barticle}[mr]
\bauthor{\bsnm{Meyn},~\bfnm{Sean~P.}\binits{S.P.}}
(\byear{2006}).
\btitle{Large deviation asymptotics and control variates for simulating large
functions}.
\bjournal{Ann. Appl. Probab.}
\bvolume{16}
\bpages{310--339}.
\bid{doi={10.1214/105051605000000737}, issn={1050-5164}, mr={2209344}}
\bptok{imsref}%
\end{barticle}
%
\endbibitem

%b15 #&#
\bibitem{Rogers00}
%
\begin{bbook}[mr]
\bauthor{\bsnm{Rogers},~\bfnm{L.~C.~G.}\binits{L.C.G.}} \AND
\bauthor{\bsnm{Williams},~\bfnm{David}\binits{D.}}
(\byear{2000}).
\btitle{Diffusions, {M}arkov Processes, and Martingales. {V}ol. 1, Foundations}.
\bseries{Cambridge Mathematical Library}.
\blocation{Cambridge}: \bpublisher{Cambridge Univ. Press}.
\bid{mr={1796539}}
\bptok{imsref}%
\end{bbook}
%
\endbibitem

\end{thebibliography}
\end{document}